\documentclass[reqno,a4paper,12pt]{amsart} 

\usepackage{amsmath,amscd,amsfonts,amssymb}
\usepackage{mathrsfs,dsfont}

\usepackage{tikz}
\usetikzlibrary{patterns}
\usepackage{float}
\usepackage{appendix}
\usepackage{caption}
\usepackage{subcaption}
\usepackage{enumerate}
\usepackage{hyperref}

\numberwithin{equation}{section}
\numberwithin{figure}{section}

\def\R{\mathbb{R}}

\def\Z{\mathbb{Z}}

\def\1{\mathds{1}}

\renewcommand\le{\leqslant}
\renewcommand\ge{\geqslant}
\renewcommand\leq{\leqslant}
\renewcommand\geq{\geqslant}

\newcommand{\diam}{\operatorname{diam}}
\newcommand{\dist}{\operatorname{dist}}

\theoremstyle{plain}
\newtheorem{thm}{Theorem}[section]
\newtheorem{lem}[thm]{Lemma}
\newtheorem{corollary}[thm]{Corollary}

\newtheorem*{claim*}{Claim}
\newtheorem*{thm*}{Theorem}

\theoremstyle{definition}
\newtheorem{definition}[thm]{Definition}
\newtheorem*{definition*}{Definition}
\newtheorem*{remarks*}{Remarks}
\newtheorem*{remark*}{Remark}

\newenvironment{enumerate-math}
{\begin{enumerate}
		\addtolength{\itemsep}{5pt}
		}
	{\end{enumerate}}

\newenvironment{enumerate-text}
{\begin{enumerate}
		\addtolength{\itemsep}{5pt}
		}
	{\end{enumerate}}

\begin{document}
	
	\title{Some constructions of restricted Kakeya sets}

    \author{Cheng Chen}
\address{Department of Mathematics \& International Center for Mathematics, Southern University of Science and Technology, Shenzhen 518055, China}
\email{charliot.cc@gmail.com}

\author{Jianbin Li}
\address{Department of Mathematics, Southern University of Science and Technology, Shenzhen 518055, China}
\email{leeanalysing@gmail.com}

\author{Longhui Li}
\address{Department of Mathematics, Southern University of Science and Technology, Shenzhen 518055, China}
\email{12231267@mail.sustech.edu.cn}

	\author{Bochen Liu}
	\address{Department of Mathematics \& International Center for Mathematics, Southern University of Science and Technology, Shenzhen 518055, China}
	\email{Bochen.Liu1989@gmail.com}

    \author{Ruichen Sun}
\address{School of Mathematics and Statistics, Northwestern Polytechnical University, Xi'an 710129, China}
\email{ruichensun@mail.nwpu.edu.cn}

    \author{Haoran Wang}
\address{Department of Mathematics, Southern University of Science and Technology, Shenzhen 518055, China}
\email{12313120@mail.sustech.edu.cn}

    \author{Saijie Zhou}
\address{Department of Mathematics, Southern University of Science and Technology, Shenzhen 518055, China}
\email{zhousaijie2003@gmail.com}

	\thanks{This project was conducted during the 2026 SUSTech Undergraduate Summer Research Program in Real Analysis. Bochen Liu was partially supported by the National Key R\&D Program of China 2024YFA1015400. The project was partially supported by SUSTech Other Teaching Funds (Y01281817), SUSTech Teaching‑Experiment Consumables and Travel Expenses (Y01281811)}
	\date{}
	
	
	\begin{abstract}
		
		In this paper, we consider Kakeya sets with the additional restriction that centers of the unit line segments belong to a given set. In particular, for every uncountable Borel set $A\subset\R^d$, $d\geq 2$, we construct a compact subset of $\R^d$ of Lebesgue measure zero that contains, in every direction, a unit line segment whose center lies in $A$. Notice that every such Kakeya set (not even necessarily compact) must have positive Lebesgue measure if $A$ is countable. So our result shows that the countability is in fact the only obstruction.
		
	\end{abstract}
	\maketitle
	
	\section{Introduction}
	
	The Kakeya problem is one of the most popular math problems recently. A set in $\R^d$ is called a Kakeya set (also called a Besicovitch set), if it contains a unit line segment in every direction. In 1928, Besicovitch \cite{Besicovitch28} constructed Kakeya sets of Lebesgue measure zero. We also refer to \cite{Fal86}\cite{SS05} for classical constructions and \cite{LW25+}\cite{FS26+} for recent constructions.
	
	Though a Kakeya set may have Lebesgue measure zero, it is still expected to be ``large". The famous Kakeya conjecture states that the Hausdorff dimension of every Kakeya set in $\R^d$ must be $d$. It has been solved in the plane by Davies \cite{Davies71} (see \cite{Wol03} for other proofs), and in $\R^3$ by Wang-Zahl \cite{WZ25+} recently.
	
	In a very recent work, Fraser and Yang \cite{FY25+} studied restricted Kakeya sets, which are Kakeya sets with the additional restriction that centers of the unit line segments belong to a given set $A\subset\R^d$. By applying Bourgain's bush argument \cite{Bou91}, they are able to obtain better estimates on dimension of restricted Kakeya sets when $A$ has smaller dimension. 
	
	Inspired by the work of Fraser-Yang and the raising of the Kakeya conjecture, we are wondering if a restricted Kakeya set can have Lebesgue measure zero, especially when $A$ is small. While a countable set $A$ is clearly insufficient, we show that the countability is in fact the only obstruction.
	
	\begin{thm}\label{main-thm-Rd}
		Suppose $A\subset\R^d$, $d\geq 2$ is an uncountable Borel set. Then there exists a compact subset of $\R^d$ of Lebesgue measure zero that contains, in every direction, a unit line segment whose center lies in $A$.
	\end{thm}
	
	It suffices to show the case $d=2$. For $d\geq 3$, there exists a two-dimensional subspace, say $\R^2\times\{0\}$, such that the projection of $A$ to this subspace is still uncountable, denoted by $A'$. We may also assume that $A\subset [0,1]^d$. If the case $d=2$ is proved, then there exists a compact $E\subset\R^2$ of Lebesgue measure zero containing unit line segments $l\in\mathcal{L}_E$ of every direction with centers $c(l)\in A'$. For every such $c(l)\in A'$, there exists some point $(c(l), c'(l))\in A$. Therefore, because of
    $$\bigcup_{l\in\mathcal{L}_E} \left(l\times B(c'(l), \frac{1}{2})\right)\subset E\times [-1, 2]^{d-2},$$
    the set $E\times [-1, 2]^{d-2}$ is our desired restricted Kakeya set in $\R^d$ of Lebesgue measure zero.
    
    We shall handle the cases whether $A$ lies in a line separately. Moreover, when $A$ lies in a line, we in fact prove a stronger result. 
    
    We remind the readers that every uncountable Borel set contains a perfect subset (see Section \ref{subset-perfect} below).
	
	\begin{thm}\label{main-thm-in-a-line}
		Suppose $0\in A\subset[0,1]$ is a perfect set. Then there exists a subset $E\subset \R^2$ of Lebesgue measure zero that contains, in every direction in $S^1\setminus\{(\pm1,0)\}$, a unit line segment centered in $A\times\{0\}$. Moreover,
        $$E\cap(\mathbb R\times\{0\})=A\times\{0\}$$
        and $E\cup ([-1,2]\times\{0\})$ is compact.
	\end{thm}
    
    If $A$ is not in a line, one cannot expect $A\subset E$ as Theorem \ref{main-thm-in-a-line} because $A$ itself may have positive Lebesgue measure. Our result for general $A$ is the following, which immediately implies Theorem \ref{main-thm-Rd} by Lemma \ref{lem:A-H perfect} below and the discussion after Theorem \ref{main-thm-Rd}.
	\begin{thm}\label{main-thm-general}
		Suppose $A\subset\R^2$ is a perfect set. Then there exists a compact subset $E\subset \R^2$ of Lebesgue measure zero that contains, in every direction, a unit line segment whose center lies in $A$.	
	\end{thm}

\section{Preliminary}
	
	In this section, we collect several standard definitions and fundamental results from geometric measure theory and point-set topology that will be used throughout the paper.
	
	\subsection{Densities and irregular sets}
	Let $\mathcal{H}^s$ denote the $s$-dimensional Hausdorff measure on $\mathbb{R}^d$, and let $\mathcal{L}^d$ denote the $d$-dimensional Lebesgue measure. For $s = 1$ in the plane $\mathbb{R}^2$, a Borel set $E \subset \mathbb{R}^2$ with $0 < \mathcal{H}^1(E) < \infty$ is called a \emph{$1$-set}.
	
	For a $1$-set $E \subset \mathbb{R}^2$ and a point $x \in \mathbb{R}^2$, the \emph{lower $1$-density} and \emph{upper $1$-density} of $E$ at $x$ are defined respectively by
	\[
	\underline{D}^1(E, x) := \liminf_{r \to 0^+} \frac{\mathcal{H}^1(E \cap B(x, r))}{2r},
	\]
	\[
	\overline{D}^1(E, x) := \limsup_{r \to 0^+} \frac{\mathcal{H}^1(E \cap B(x, r))}{2r},
	\]
	where $B(x, r)$ denotes the closed ball centered at $x$ with radius $r > 0$.
	
	\begin{definition}
		A $1$-set $E \subset \R^2$ is called  \emph{irregular} if
		\[
		\underline{D}^1(E, x) < 1 \quad \text{for } \mathcal{H}^1\text{-almost every } x \in E.
		\]
	\end{definition}

	For each $t \in \mathbb{R}$, we define the projection mapping $\pi_t : \mathbb{R}^2 \to \mathbb{R}$ by 
	\[
	\pi_t(b, a) := bt + a,
	\]
	and for the vertical direction, let $\pi_\infty(b, a) := b$. Note that $\pi_t(b, a)$ coincides with the standard orthogonal projection of the point $(b, a)$ onto the line in direction $(t, 1)$, up to a non-zero scaling factor $\sqrt{1+t^2}$.  Under this definition, the classical Besicovitch projection theorem for irregular sets is stated as follows:
	
	\begin{thm}[see, e.g. {\cite[Theorem 6.13]{Fal86}}\label{thm:besicovitch}]
		Let $E \subset \mathbb{R}^2$ be an irregular $1$-set. Then the projection $\pi_t(E)$ has $1$-dimensional Lebesgue measure zero for almost every $t \in \mathbb{R}$.
	\end{thm}
	
	\subsection{Perfect sets}\label{subset-perfect}
	A subset $A \subset \mathbb{R}^d$ is called a \emph{perfect set} if it is closed and contains no isolated points.
	
	\begin{lem}[Alexandrov--Hausdorff, see, e.g. Theorem 13.6 in \cite{Kechris95}]\label{lem:A-H perfect}
		Every uncountable Borel set in $\R^d$ contains a compact subset homeomorphic to the ternary Cantor set, i.e., a non-empty, bounded, totally disconnected perfect set.
	\end{lem}

    Every perfect set is uncountable. In fact, if $A$ is perfect and $x\in A$, then every neighborhood of $x$ contains uncountably many points in $A$, for otherwise, if $A=\{x_i\}$, then one can construct a decreasing sequence of closed balls $\{B_i\}$ with $B_i\cap A\neq \emptyset$ and $x_i\notin B_i$. So Lemma \ref{lem:A-H perfect} has the following corollary.
	
	\begin{corollary} \label{cor:CB_extract}
		Let $A \subset \R^n$ be a perfect set. If $U$ is a non-empty relatively open subset of $A$, then $U$ contains a non-empty perfect subset.
	\end{corollary}

	\section{Proof of Theorem \ref{main-thm-in-a-line}}
	We use the point-line duality to produce Kakeya sets with zero Lebesgue measure. This idea has been used on classical constructions of Kakeya sets. We refer to Chapter 7 in \cite{Fal86} as a reference.
    
    For $(b,a)\in \R^2$, we define the affine line $l_{b,a}$ by
	\begin{align*}
		l_{b,a}=\{(a,0)+t(b,1):t\in\R\}.
	\end{align*}
	
	We will construct a compact set $D\subset \R^2$ satisfying
	\begin{enumerate}
		\item the vertical projection $\pi_\infty(D) = [0,1]$ ;
		\item the horizontal projection $\pi_0(D) = A$ ;
		\item $D$ is an irregular $1$-set,
	\end{enumerate} 
and take
	\[
	E_0:= \left\{ (a, 0) + t(b, 1) : (b,a) \in D, \; t \in \left[-\frac{1}{2}, \frac{1}{2}\right] \right\}.
	\]
Clearly $E_0$ contains unit line segments in all directions in $\{\frac{(b,1)}{\sqrt{b^2+1}}, b\in[0,1]\}$ centered in $A\times\{0\}$ and $E_0\cap (\R\times\{0\})=A\times\{0\}$. To compute the Lebesgue measure of $E_0$, consider each horizontal slice of $E_0$ at $y \in \R$:
	\[
	E_0^y := \{ by+a: (b, a) \in D \}=\pi_y(D).
	\]
	As $D$ is an irregular $1$-set, Theorem~\ref{thm:besicovitch} implies that $\mathcal{L}^1(E_0^y) = 0$ for almost every $y \in \mathbb{R}$. Therefore, by Fubini's theorem,
	\[
	\mathcal{L}^2(E_0) = \int_{-\infty}^\infty \mathcal{L}^1(E_0^y) \, dy = 0.
	\]
	
	Finally, by considering $D+(n,0), n\in\Z$ and letting
    \[
	E:= \bigcup_{n\in\Z}E_n:=\bigcup_{n\in\Z}\left\{ (a, 0) + t(b, 1) : (b,a) \in D+(n,0), \; t \in \left[-\frac{1}{2}, \frac{1}{2}\right] \right\},
	\]
    we obtain a planar set of Lebesgue measure zero that contains, in every direction in $S^1 \setminus \{(\pm 1, 0)\}$, a unit line segment centered in $A\times \{0\}$. Also 
    $$E \cap (\mathbb{R} \times \{0\}) = \bigcup_{n\in\Z} E_n \cap (\mathbb{R} \times \{0\})=A\times \{0\}. $$
    To see $E\cup ([-1,2]\times\{0\})$ is compact, notice that, for every convergent sequence in $E\cup ([-1,2]\times\{0\})$, either there exists a subsequence contained in one of $\{[-1,2]\times\{0\}, E_n, n\in \Z\}$, or there exists a subsequence $p_i\in E_{n_j}$ with $|n_j|\rightarrow\infty$. In either case, the limit lies in $E\cup ([-1,2]\times\{0\})$.

	Now we turn to the construction of $D$. We construct $D$ as the intersection of a nested sequence of compact sets $D_0 \supset D_1 \supset D_2 \supset \cdots$, where $D_0 = [0,1]^2$, $D_k$ is a union of disjoint closed squares $Q$ of side length $N_k^{-1}$, $N_k\in\Z_+$, satisfying $\pi_\infty(D_k)=[0,1]$ and the interior of $\pi_0(Q)\cap A$ for each $Q$ is nonempty.

Assume $D_k$ is constructed. At step $k+1$, let $Q = I_Q \times J_Q \subset D_k$ be a $N_k^{-1}$ square. We partition $Q$ into subsquares of side length $N_{k+1}^{-1} = (N_{k+1}' \cdot N_k)^{-1}$. Because $A$ is perfect and the interior of $\pi_0(Q)\cap A$ is nonempty, for sufficiently large $N_{k+1}$, there are $4\leq M_{Q} \leq N_{k+1}'$ many $N_{k+1}^{-1}$-rows in $Q$ whose interiors intersect $\{0\}\times A$. We call them ``valid" and denote their horizontal projections by $\{J_Q^{(j)}\}_{j=1}^{M_{Q}}$ from bottom to top. Also denote vertical projections of all $N_{k+1}^{-1}$-columns in $Q$ as $\{I_Q^{(j')}\}_{j'=1}^{N_{k+1}'}$ from left to right. As there are only finitely many cubes in each step, we may assume $N_{k+1}$ is universal over all $N_k^{-1}$-cubes in $D_k$.

For convenience, we also use $J_Q^{(j)}$ to denote the row $\R\times J_Q^{(j)}$, and use $I_Q^{(j')}$ to denote the column $I_Q^{(j')}\times\R$.
    
	Now we make selections over $N_{k+1}^{-1}$-squares $\{I_Q^{(j')}\times J_Q^{(j)}: 1\leq j\leq M_Q, 1\leq j'\leq N_{k+1}'\}$.
    
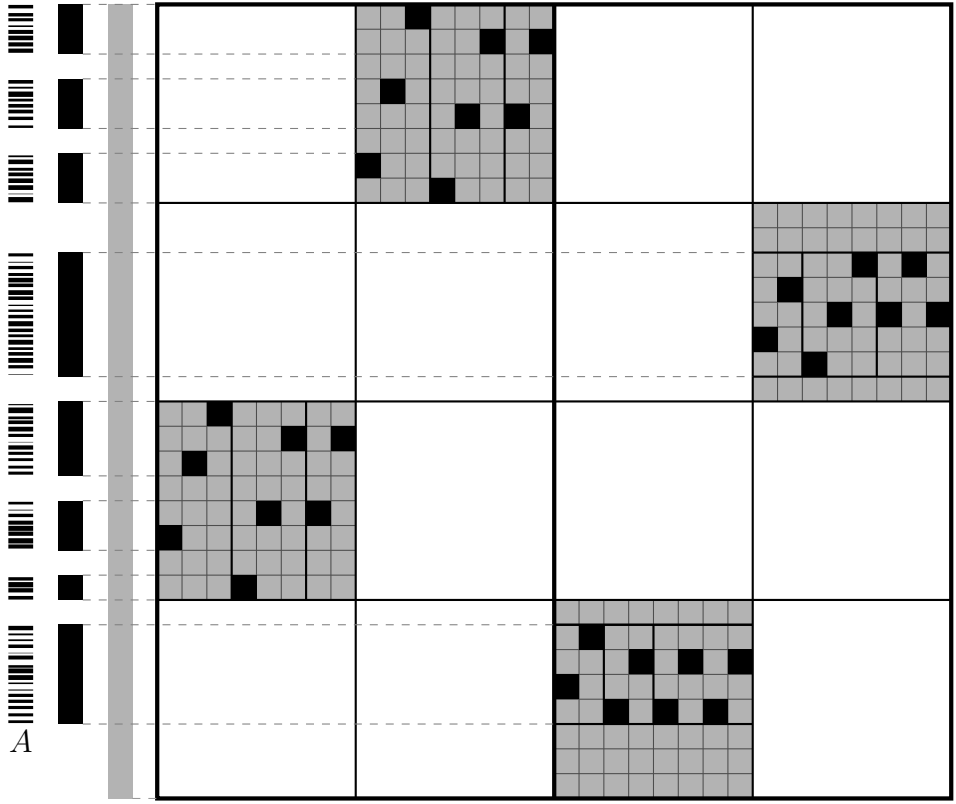
\begin{figure}[htbp]
		\centering
		\begin{tikzpicture}[scale=0.75]
			
			\def\l{3.5/8}	
			
			\foreach \x/\y in {-7/3.5,-3.5/10.5,0/0,3.5/7}
			{
				\fill[black!30] (\x,\y) rectangle ++(3.5,3.5);
				\foreach \n in {1,...,7}
				{
					\draw[black!70] ({\x+\n*\l},\y)--++(0,3.5);
					\draw[black!70] (\x,{\y+\n*\l})--++(3.5,0);
				}
			}
			
			\fill ({-7+3*\l},3.5) rectangle ++(\l,\l);
			\fill ({-7},{3.5+2*\l}) rectangle ++(\l,\l);
			\fill ({-7+4*\l},{3.5+3*\l}) rectangle ++(\l,\l);
			\fill ({-7+\l},{3.5+5*\l}) rectangle ++(\l,\l);
			\fill ({-7+5*\l},{3.5+6*\l}) rectangle ++(\l,\l);
			\fill ({-7+2*\l},{3.5+7*\l}) rectangle ++(\l,\l);
			\fill ({-7+6*\l},{3.5+3*\l}) rectangle ++(\l,\l);
			\fill ({-7+7*\l},{3.5+6*\l}) rectangle ++(\l,\l);
			\foreach \j in {1,2}
			\draw[thick] ({-7+3*\j*\l},3.5)--++(0,3.5);

			\foreach \j in {0,1}
			\fill ({3.5+\j*\l},{7+2*\l+2*\j*\l}) rectangle ++(\l,\l);
			\foreach \j in {0,1,2}
			\fill ({3.5+(2+\j)*\l},{7+\l+2*\j*\l}) rectangle ++(\l,\l);
			\foreach \j in {1,3}
			\fill ({7-\j*\l},{7+3*\l}) rectangle ++(\l,\l);
			\fill ({7-2*\l},{7+5*\l}) rectangle ++(\l,\l);
			\draw[thick] (3.5,{7+\l})--++(3.5,0);
			\draw[thick] (3.5,{7+6*\l})--++(3.5,0);
			\draw[thick] ({3.5+2*\l},{7+\l})--++(0,5*\l);
			\draw[thick] ({3.5+5*\l},{7+\l})--++(0,5*\l);
			
			\fill (0,{3.5-4*\l}) rectangle ++(\l,\l);
			\fill (\l,{3.5-2*\l}) rectangle ++(\l,\l);
			\foreach \n in {2,4,6}
			{
				\fill ({\n*\l},{3.5-5*\l}) rectangle ++(\l,\l);
				\fill ({(\n+1)*\l},{3.5-3*\l}) rectangle ++(\l,\l);
			}
			\draw[thick] (0,{3.5-\l})--++(3.5,0);
			\foreach \j in {1,2}
			\draw[thick] (2*\j*\l,3.5-\l)--++(0,-4*\l);
			\draw[thick] (0,{3.5-\l})--++(3.5,0);
			\draw[thick] (0,{3.5-5*\l})--++(3.5,0);
			
			\foreach \k in {0,1}
			\foreach \n in {0,1,2}
			\fill ({-3.5+\n*\l+3*\k*\l},{10.5+(3*\n+1-\k)*\l}) rectangle ++(\l,\l);
			\fill (-2*\l,{10.5+3*\l}) rectangle ++(\l,\l);
			\fill (-\l,14-2*\l) rectangle ++(\l,\l);
			\foreach \j in {1,2}
			\draw[thick] ({-3.5+3*\j*\l},10.5)--++(0,3.5);
			
			\draw[ultra thick] (-7,0) rectangle (7,14);
			\draw[ultra thick] (0,0)--(0,14);
			
			\foreach \n in {-3.5,0,3.5}
			{
				\draw[thick] (\n,0)--++(0,14);
				\draw[thick] (-7,\n+7)--++(14,0);
			}
			
			\fill[black!30] ({-7-2*\l},0) rectangle++(\l,14);
			\draw[black!40, dashed] ({-7-\l},0)--(-7,0);
			
			\fill({-7-4*\l},3*\l) rectangle ++(\l,4*\l);
			\fill({-7-4*\l},3.5) rectangle ++(\l,\l);
			\fill({-7-4*\l},{3.5+2*\l}) rectangle ++(\l,2*\l);
			\fill({-7-4*\l},{3.5+5*\l}) rectangle ++(\l,3*\l);
			\fill ({-7-4*\l},{7+\l}) rectangle ++(\l,5*\l);
			\foreach \j in {0,3,6}
			\fill ({-7-4*\l},{10.5+\j*\l}) rectangle ++(\l,2*\l);
			
			\foreach \j in {3,7}
			\draw[black!50, dashed] ({-7-3*\l},\j*\l)--(0,\j*\l);
			\foreach \j in {0,1,2,4,5,8,16,24}
			\draw[black!50, dashed]
			({-7-3*\l},{3.5+\j*\l})--(-7,3.5+\j*\l);
			\foreach \j in {1,6}
			\draw[black!50, dashed] ({-7-3*\l},{7+\j*\l})--(3.5,{7+\j*\l});
			\foreach \j in {2,3,5,6}
			\draw[black!50, dashed] ({-7-3*\l},{10.5+\j*\l})--(-3.5,{10.5+\j*\l});
			
			\fill \foreach \k in {3,4,5,6, 8, 10,11, 13,14,15, 17,18,19,20,21, 24,25, 27,28, 30,31} {
				({-7-6*\l}, {(\k + 0.01 + 0.08*abs(sin(\k*113)))*\l}) rectangle ({-7-5*\l}, {(\k + 0.10 + 0.14*abs(cos(\k*37)))*\l})
				({-7-6*\l}, {(\k + 0.26 + 0.09*abs(sin(\k*71)))*\l}) rectangle ({-7-5*\l}, {(\k + 0.36 + 0.13*abs(cos(\k*199)))*\l})
				({-7-6*\l}, {(\k + 0.50 + 0.11*abs(sin(\k*23)))*\l}) rectangle ({-7-5*\l}, {(\k + 0.62 + 0.14*abs(cos(\k*107)))*\l})
				({-7-6*\l}, {(\k + 0.77 + 0.08*abs(sin(\k*89)))*\l}) rectangle ({-7-5*\l}, {(\k + 0.86 + 0.13*abs(cos(\k*53)))*\l})
			};
			\node at ({-7-5.5*\l},{-0.3+3*\l}){$A$};
			
			\node at ({-7-1.5*\l},-0.3) {};
			
		\end{tikzpicture}
		
		\caption{{Two generations inside the parent square $C_0$}
		}
		\label{fig:construction}
	\end{figure}

    In columns $I_Q^{(j')}, j' \in \{1, \dots, N_{k+1}'\}$, we select exactly one $N_{k+1}^{-1}$-square in each column using a parity-alternating pattern. The alternating pattern ensures that adjacent columns are separated by vertical gaps, so that when computing the density at small scales a ball of radius $r_k$ can intersect only one parity family of selected squares.  More precisely, from left to right:
	\begin{itemize}
		\item For the first $\lfloor M_{Q}/2 \rfloor$ columns, select squares in the even rows $J_Q^{(2)}, J_Q^{(4)}, \dots$
		\item For the next $\lceil M_{Q}/2 \rceil$ columns, select squares in the odd rows $J_Q^{(1)}, J_Q^{(3)}, \dots$
		\item Since $M_{Q}\leq N_{k+1}'$, all valid rows will be selected. For remaining columns, alternate between the last two odd valid rows. 
	\end{itemize}
    See Figure \ref{fig:construction} below for an illustration of this selection.

	Let $D_{k+1}$ be the union of all selected $N_{k+1}^{-1}$-subsquares across all parent $N_k^{-1}$-squares $Q \subset D_k$. We define the limit set by $D := \bigcap_{k=0}^\infty D_k$.
	
	By standard compactness and nested intersection arguments:
	\begin{enumerate}
		\item Since one square is chosen in each column in every step, the vertical projection keeps surjective: $\pi_\infty(D) = \pi_\infty(D_k)=[0,1]$.
		\item Since only rows intersecting $A$ are retained and every valid row is selected, the horizontal projection of each $D_k$ satisfies $A\subset \pi_0(D_k)\subset \mathcal{N}_{N_k^{-1}}(A)$, which implies $\pi_0(D) = A$ as $A$ is compact.
	\end{enumerate}
	
	It remains to prove $D$ is an irregular $1$-set. Since $\pi_\infty(D)=[0,1]$, the lower bound on $\mathcal H^1(D)$ follows immediately from the fact that $\pi_\infty$ is $1$-Lipschitz:
	\[
	\mathcal H^1(D)\ge \mathcal H^1\big(\pi_\infty(D)\big)>0.
	\]
	To prove the upper bound, since there is exactly one $N_{k+1}^{-1}$-square is selected in each $N_{k+1}^{-1}$-column, the number of $N_{k+1}^{-1}$-squares in $D_{k+1}$ is exactly $N_{k+1}$. Therefore,
	\[
	\mathcal{H}^1(D)\le\liminf_{k\to\infty}\sqrt{2}\, N_{k+1}^{-1}\cdot N_{k+1}=\sqrt{2}.
	\] 
	Hence $D$ is a $1$-set.
	
	To verify irregularity, consider an arbitrary $p \in D$ and $r_k := N_k^{-1}$. For every $l>k$, let $\mathcal Q_l$ denote the family of all
	$N_l^{-1}$-squares of $D_l$ intersecting $B(p,r_k)$. Since distinct parent squares are separated by at least their side length, the ball $B(p, r_k)$ intersects at most one parent square of generation $k$. As a consequence, by column-counting, the number of subsquares of generation $l > k$ intersecting $B(p, r_k)$ satisfies $\#\mathcal{Q}_l \le r_k / N_l^{-1}+2$. 
	Since $D\subset D_l$, the family $\mathcal Q_l$ covers $D\cap B(p,r_k)$. So we have the measure estimate
	\[
	\mathcal{H}^1(D \cap B(p, r_k)) \le	\liminf_{l\to\infty}
	\sqrt2\,N_l^{-1}\#\mathcal Q_l= \sqrt{2} r_k,
	\]
	and then
	\[ \underline{D}^1(D, p)\le \liminf_{k\to\infty}
	\frac{\mathcal H^1(D\cap B(p,r_k))}{2r_k} \le \frac{\sqrt{2}}{2} < 1.\]
	Hence, $D$ is an irregular $1$-set.

	\section{Proof of Theorem \ref{main-thm-general}}
	Now we deal with $A\subset\R^2$ that is not necessarily contained in a line. By Corollary \ref{cor:CB_extract}, we may assume $A$ is compact. We may also assume that the projection of $A$ to every direction is a perfect set. This is because, if the projection of $A$ in the direction of $e$ contains an isolated point, then there exists a line $l$ perpendicular to $e$ such that a neighborhood of $A\cap l$ contains no points from $A\backslash l$. Therefore $A\cap l$ is a perfect set, and one can take the construction from the previous section.
	
	We shall construct a compact set $D\subset \R^2$ satisfying
	\begin{enumerate}
		\item for every $b\in[0,1]$, there exists $(x, y)\in A$ such that $(b,x-yb)\in D$,
		\item $D$ is an irregular $1$-set.
	\end{enumerate}
	Then we take 
	\[
		E_0:=\{(x,y)+t(b,1): (x,y)\in A, b\in[0,1], (b,x-yb)\in D, t\in[-\frac{1}{2}, \frac{1}{2}]\},
	\]
    which is a compact set that contains, in every direction in $\{\frac{(b,1)}{\sqrt{b^2+1}}, b\in[0,1]\}$, a unit line segment centered in $A$ . For its Lebesgue measure, notice that each horizontal $h$-slice
	$$E_0^h:=\{x+(h-y)b: (b,x-yb)\in D\}=\{bh+(x-yb): (b,x-yb)\in D\}$$
	is contained in $\pi_{h}(D)$ that has Lebesgue measure zero for almost all $h$ by Theorem \ref{thm:besicovitch}. Hence $E_0$ has Lebesgue measure zero by the Fubini theorem. 
	
	By repeating Corollary \ref{cor:CB_extract}, one can find a sequence of compact perfect sets $A\supset A_1\supset A_2\supset\cdots$ such that $\diam (A_n)\rightarrow 0$ and $\cap A_n=\{(x_0,y_0)\}\in A$. For each $A_n$ in place of $A$, there is an associated $D_{A_n}$ as above. Then, similar to the previous section, by taking $D'_{A_n}:=D_{A_n}+(n,0)$,
    $$
		E_n:=\{(x,y)+t(b,1): (x,y)\in A_n, b\in[n,n+1], (b,x-yb)\in D_{A_n}+(n,0), t\in[-\frac{1}{2}, \frac{1}{2}]\},$$
        and finally
        $$E:=\bigcup_{n\in\Z} E_n \cup\{(x_0, y_0)+t(1,0): t\in[-\frac{1}{2}, \frac{1}{2}]\},$$
   we obtain a compact subset $E\subset \R^2$ of Lebesgue measure zero that contains, in every direction, a unit line segment centered in $A$. Here $E$ is compact because $E_n\rightarrow \{(x_0, y_0)+t(1,0): t\in[-\frac{1}{2}, \frac{1}{2}]\}$ under the Hausdorff metric.
	
	It remains to construct $D$. Similar to the previous section, $D$ is the intersection of a nested sequence of compact sets $D_0\supset D_1\supset\cdots$. We shall see that, because $y$ may not be zero for $(x,y)\in A$, the condition $(b, x-yb)\in D$ makes the construction much more complicated than the previous section.

	Suppose $(\frac{1}{4},\frac{1}{8})\in A$, then the set $A\cap \{(x,y):\frac{1}{2}>x>y>0\}$ is a non-empty relatively open subset of $A$, which by Corollary~\ref{cor:CB_extract}, contains a perfect subset. So we may assume $$A\subset \{(x,y):\frac{1}{2}>x>y>0\}.$$

    Let $D_0=[0,1]^2$. More generally, $D_k $ consists of $ M_k $ pairwise disjoint closed squares
	\[
	Q_{k,j} = I_{k,j} \times J_{k,j},\qquad j=1,\ldots,M_k,
	\]
	 with $ |I_{k,j}| = |J_{k,j}| = l_{k,j} $, and the vertical projection satisfies $\pi_{\infty}(D_k)=[0,1]$.  Also the following property on $Q_{k,j}$ will be preserved in every step: for every $b\in I_{k,j}$, there exists some $(x,y)\in A$ such that $x-yb\in \operatorname{Int}(J_{k,j})$. In particular, the set
     $$A_{k,j}:=\{(x,y)\in A: x-yb\in \operatorname{Int}(J_{k,j}) \text{ for some } b\in I_{k,j}\}$$
     is a nonempty relatively open subset of $A$, thus must be uncountable and contains a perfect set by Corollary \ref{cor:CB_extract}. So we may assume $A_{k,j}$ does not lie in a line, otherwise one can simply take the construction from Theorem \ref{main-thm-in-a-line}.

     Now we construct $D_{k+1}$.

	Note that for a fixed point $(x,y)\in A$, as $b$ varies, 
	\begin{equation}\label{eq-w}
	w_{(x,y)}(b) := x - yb
	\end{equation}
	describes a straight line segment in the parameter space $(b,w)$. Therefore, as $ A_{k,j} $ is uncountable and not contained in a line, for any $ b \in I_{k,j}$ there exist two distinct points $ (x_b^+, y_b^+), (x_b^-, y_b^-) \in A_{k,j} $ such that
	$$
	w_{(x_b^+, y_b^+)} (b) > w_{(x_b^-, y_b^-)} (b), \text{ both in } \operatorname{Int}(J_{k,j}).
	$$
	By the continuity of $w$, they remain strictly separated and away from endpoints of $J_{k,j}$ on an open neighborhood $U_b$ of $b$. The collection $\{U_b\}_{b \in I_{k,j}}$ forms an open cover of the compact interval $I_{k,j}$.
	
	By compactness, there exists a finite subcover $ \{U_{b_1}, \ldots, U_{b_N}\} $ of $I_{k,j}$, and a partition
    $$I_{k,j}= \bigcup_{m=1}^{S_j}[c_{m-1}, c_m]:=\bigcup_{m=1}^{S_j}H_{j,m},\ 
	c_0 < c_1 < \cdots < c_{S_j},$$
such that each closed interval $H_{j,m}$ is contained in some $U_{b_n}$. So on each $ H_{j,m} $, one can find a pair of points $ (x_{j,m}^\pm, y_{j,m}^\pm)\in A_{k,j} $ such that
	\[
	w_{(x_{j,m}^+, y_{j,m}^+)} (b)>w_{(x_{j,m}^-, y_{j,m}^-)} (b), \text{ both in } \operatorname{Int}(J_{k,j}), \ \forall\, b\in H_{j,m}.
	\]
	Thus, by the compactness of closed intervals and the finiteness of $j, m$, one can conclude that
	\[
	\Delta_k := \min_{1 \le j \le M_k} \min_{1 \le m \le S_j} \min_{b\in H_{j,m}}|w_{(x_{j,m}^+, y_{j,m}^+)} (b)-w_{(x_{j,m}^-, y_{j,m}^-)} (b)| > 0
	\]
	and all $w_{(x_{j,m}^\pm, y_{j,m}^\pm)} (b)$ are away from endpoints of $J_{k,j}$:
    $$ \min_{1 \le j \le M_k} \min_{1 \le m \le S_j} \min_{b\in H_{j,m}}\dist (w_{(x_{j,m}^\pm, y_{j,m}^\pm)} (b), \partial (J_{k,j})) > 0.$$

	Having established the global gap $\Delta_k$, we divide the $H_{j,m}$ into $N_{k+1}$ intervals $I_{j,m}^{(i)}, i = 1, \dots, N_{k+1}$, of equal length $l_{k+1, j, m} = |H_{j,m}| / N_{k+1}$. By choosing a sufficiently large integer $N_{k+1}$, one can make the maximum child width $l_{k+1}^{\max}:=\max\limits_{j,m}\{l_{k+1,j,m}\}$ satisfy:
	\begin{equation} \label{eq:scale2}
		\left( \sum_{j=1}^{M_k} S_j \right) l_{k+1}^{\max} \le \frac{0.49 \Delta_k}{100}.
	\end{equation}
	Let $b_{j,m}^{(i)}$ denote the midpoint of $I_{j,m}^{(i)}$, and then
    $$J_{j,m}^{(i)} = \begin{cases}
        [ w_{(x_{j,m}^+, y_{j,m}^+)}(b_{j,m}^{(i)}) - \frac{l_{k+1, j, m}}{2}, w_{(x_{j,m}^+, y_{j,m}^+)}(b_{j,m}^{(i)}) + \frac{l_{k+1, j, m}}{2} ], & \text{ if $i$ is odd },\\ [ w_{(x_{j,m}^-, y_{j,m}^-)}(b_{j,m}^{(i)}) - \frac{l_{k+1, j, m}}{2}, w_{(x_{j,m}^-, y_{j,m}^-)}(b_{j,m}^{(i)}) + \frac{l_{k+1, j, m}}{2} ], & \text{ if $i$ is even}
    \end{cases}.$$
    Since all $w_{(x_{j,m}^\pm, y_{j,m}^\pm)} (b_{j,m}^{(i)})$ are away from endpoints of $J_{k,j}$, one can choose $N_{k+1}$ large such that all resulting child blocks $R_{j,m}^{(i)} = I_{j,m}^{(i)} \times J_{j,m}^{(i)}$ are contained in $Q_{k,j}$.

\begin{figure}[htbp]
		\centering
		\begin{tikzpicture}[scale=1.1, every node/.style={font=\small}]
			\def\side{10}               
			
			\pgfmathsetmacro{\wone}{10*2/16}
			\pgfmathsetmacro{\wtwo}{10*3/16}
			\pgfmathsetmacro{\wthree}{10*2/16}
			\pgfmathsetmacro{\wfour}{10*4/16}
			\pgfmathsetmacro{\wfive}{10*5/16}
			
			\pgfmathsetmacro{\xone}{\wone}
			\pgfmathsetmacro{\xtwo}{\wone + \wtwo}
			\pgfmathsetmacro{\xthree}{\wone + \wtwo + \wthree}
			\pgfmathsetmacro{\xfour}{\wone + \wtwo + \wthree + \wfour}
			
			\pgfmathsetmacro{\midone}{\wone/2}
			\pgfmathsetmacro{\midtwo}{\wone + \wtwo/2}
			\pgfmathsetmacro{\midthree}{\wone + \wtwo + \wthree/2}
			\pgfmathsetmacro{\midfour}{\wone + \wtwo + \wthree + \wfour/2}
			\pgfmathsetmacro{\midfive}{\wone + \wtwo + \wthree + \wfour + \wfive/2}
			
			\def\colsPerStrip{4}
			
			\draw[ultra thick] (0,0) rectangle (\side,\side);
			
			\draw[thick] (\xone,0) -- (\xone,\side);
			\draw[thick] (\xtwo,0) -- (\xtwo,\side);
			\draw[thick] (\xthree,0) -- (\xthree,\side);
			\draw[thick] (\xfour,0) -- (\xfour,\side);
			
			\pgfmathsetmacro{\stepone}{\wone/\colsPerStrip}
			\foreach \i in {1,...,3} {
				\pgfmathsetmacro{\xx}{\i*\stepone}
				\draw[gray!60] (\xx,0) -- (\xx,\side);
			}
			\pgfmathsetmacro{\steptwo}{\wtwo/\colsPerStrip}
			\foreach \i in {1,...,3} {
				\pgfmathsetmacro{\xx}{\xone + \i*\steptwo}
				\draw[gray!60] (\xx,0) -- (\xx,\side);
			}
			\pgfmathsetmacro{\stepthree}{\wthree/\colsPerStrip}
			\foreach \i in {1,...,3} {
				\pgfmathsetmacro{\xx}{\xtwo + \i*\stepthree}
				\draw[gray!60] (\xx,0) -- (\xx,\side);
			}
			\pgfmathsetmacro{\stepfour}{\wfour/\colsPerStrip}
			\foreach \i in {1,...,3} {
				\pgfmathsetmacro{\xx}{\xthree + \i*\stepfour}
				\draw[gray!60] (\xx,0) -- (\xx,\side);
			}
			\pgfmathsetmacro{\stepfive}{\wfive/\colsPerStrip}
			\foreach \i in {1,...,3} {
				\pgfmathsetmacro{\xx}{\xfour + \i*\stepfive}
				\draw[gray!60] (\xx,0) -- (\xx,\side);
			}
			
			\def\startLowA{2.0}
			\def\slopeLowA{0.3}
			\def\startHighA{6.0}
			\def\slopeHighA{-0.4}
			
			\def\startLowB{5.0}
			\def\slopeLowB{-0.35}
			\def\startHighB{7.0}
			\def\slopeHighB{0.3}
			
			\def\startLowC{1.5}
			\def\slopeLowC{0.25}
			\def\startHighC{4.3}
			\def\slopeHighC{0.35}
			
			\def\startLowD{3.0}
			\def\slopeLowD{0.3}
			\def\startHighD{7.0}
			\def\slopeHighD{0.4}
			
			\def\startLowE{6.0}
			\def\slopeLowE{-0.4}
			\def\startHighE{8.0}
			\def\slopeHighE{-0.3}
			
			\draw[black, thick] (0, \startLowA) -- (\xone, \startLowA + \slopeLowA*\xone);
			\draw[black, thick] (0, \startHighA) -- (\xone, \startHighA + \slopeHighA*\xone);
			
			\pgfmathsetmacro{\tw}{\xtwo - \xone}
			\draw[black, thick] (\xone, \startLowB) -- (\xtwo, \startLowB + \slopeLowB*\tw);
			\draw[black, thick] (\xone, \startHighB) -- (\xtwo, \startHighB + \slopeHighB*\tw);
			
			\pgfmathsetmacro{\thw}{\xthree - \xtwo}
			\draw[black, thick] (\xtwo, \startLowC) -- (\xthree, \startLowC + \slopeLowC*\thw);
			\draw[black, thick] (\xtwo, \startHighC) -- (\xthree, \startHighC + \slopeHighC*\thw);
			
			\pgfmathsetmacro{\fw}{\xfour - \xthree}
			\draw[black, thick] (\xthree, \startLowD) -- (\xfour, \startLowD + \slopeLowD*\fw);
			\draw[black, thick] (\xthree, \startHighD) -- (\xfour, \startHighD + \slopeHighD*\fw);
			
			\pgfmathsetmacro{\fvw}{\side - \xfour}
			\draw[black, thick] (\xfour, \startLowE) -- (\side, \startLowE + \slopeLowE*\fvw);
			\draw[black, thick] (\xfour, \startHighE) -- (\side, \startHighE + \slopeHighE*\fvw);
			
			
			\foreach \col in {0,...,3} {
				\pgfmathsetmacro{\xcenter}{(\col + 0.5)*\stepone}
				\pgfmathsetmacro{\yLow}{\startLowA + \slopeLowA*\xcenter}
				\pgfmathsetmacro{\yHigh}{\startHighA + \slopeHighA*\xcenter}
				\ifodd\col
				\pgfmathsetmacro{\ycenter}{\yLow}
				\else
				\pgfmathsetmacro{\ycenter}{\yHigh}
				\fi
				\pgfmathsetmacro{\xstart}{\xcenter - \stepone/2}
				\pgfmathsetmacro{\ystart}{\ycenter - \stepone/2}
				\fill[gray!60] (\xstart, \ystart) rectangle ++(\stepone,\stepone);
				\draw[black] (\xstart, \ystart) rectangle ++(\stepone,\stepone);
			}
			
			\foreach \col in {0,...,3} {
				\pgfmathsetmacro{\xcenter}{\xone + (\col + 0.5)*\steptwo}
				\pgfmathsetmacro{\dx}{\xcenter - \xone}
				\pgfmathsetmacro{\yLow}{\startLowB + \slopeLowB*\dx}
				\pgfmathsetmacro{\yHigh}{\startHighB + \slopeHighB*\dx}
				\ifodd\col
				\pgfmathsetmacro{\ycenter}{\yLow}
				\else
				\pgfmathsetmacro{\ycenter}{\yHigh}
				\fi
				\pgfmathsetmacro{\xstart}{\xcenter - \steptwo/2}
				\pgfmathsetmacro{\ystart}{\ycenter - \steptwo/2}
				\fill[gray!60] (\xstart, \ystart) rectangle ++(\steptwo,\steptwo);
				\draw[black] (\xstart, \ystart) rectangle ++(\steptwo,\steptwo);
			}
			
			\foreach \col in {0,...,3} {
				\pgfmathsetmacro{\xcenter}{\xtwo + (\col + 0.5)*\stepthree}
				\pgfmathsetmacro{\dx}{\xcenter - \xtwo}
				\pgfmathsetmacro{\yLow}{\startLowC + \slopeLowC*\dx}
				\pgfmathsetmacro{\yHigh}{\startHighC + \slopeHighC*\dx}
				\ifodd\col
				\pgfmathsetmacro{\ycenter}{\yLow}
				\else
				\pgfmathsetmacro{\ycenter}{\yHigh}
				\fi
				\pgfmathsetmacro{\xstart}{\xcenter - \stepthree/2}
				\pgfmathsetmacro{\ystart}{\ycenter - \stepthree/2}
				\fill[gray!60] (\xstart, \ystart) rectangle ++(\stepthree,\stepthree);
				\draw[black] (\xstart, \ystart) rectangle ++(\stepthree,\stepthree);
			}
			
			\foreach \col in {0,...,3} {
				\pgfmathsetmacro{\xcenter}{\xthree + (\col + 0.5)*\stepfour}
				\pgfmathsetmacro{\dx}{\xcenter - \xthree}
				\pgfmathsetmacro{\yLow}{\startLowD + \slopeLowD*\dx}
				\pgfmathsetmacro{\yHigh}{\startHighD + \slopeHighD*\dx}
				\ifodd\col
				\pgfmathsetmacro{\ycenter}{\yLow}
				\else
				\pgfmathsetmacro{\ycenter}{\yHigh}
				\fi
				\pgfmathsetmacro{\xstart}{\xcenter - \stepfour/2}
				\pgfmathsetmacro{\ystart}{\ycenter - \stepfour/2}
				\fill[gray!60] (\xstart, \ystart) rectangle ++(\stepfour,\stepfour);
				\draw[black] (\xstart, \ystart) rectangle ++(\stepfour,\stepfour);
			}
			
			\foreach \col in {0,...,3} {
				\pgfmathsetmacro{\xcenter}{\xfour + (\col + 0.5)*\stepfive}
				\pgfmathsetmacro{\dx}{\xcenter - \xfour}
				\pgfmathsetmacro{\yLow}{\startLowE + \slopeLowE*\dx}
				\pgfmathsetmacro{\yHigh}{\startHighE + \slopeHighE*\dx}
				\ifodd\col
				\pgfmathsetmacro{\ycenter}{\yLow}
				\else
				\pgfmathsetmacro{\ycenter}{\yHigh}
				\fi
				\pgfmathsetmacro{\xstart}{\xcenter - \stepfive/2}
				\pgfmathsetmacro{\ystart}{\ycenter - \stepfive/2}
				\fill[gray!60] (\xstart, \ystart) rectangle ++(\stepfive,\stepfive);
				\draw[black] (\xstart, \ystart) rectangle ++(\stepfive,\stepfive);
			}
			
			
			\node[below] at (\side/2, -0.6) {$Q_{k,j}$};
			\node[below, font=\scriptsize] at (\midone, -0.4) {$H_{j,1}$};
			\node[below, font=\scriptsize] at (\midtwo, -0.4) {$H_{j,2}$};
			\node[below, font=\scriptsize] at (\midthree, -0.4) {$H_{j,3}$};
			\node[below, font=\scriptsize] at (\midfour, -0.4) {$H_{j,4}$};
			\node[below, font=\scriptsize] at (\midfive, -0.4) {$H_{j,5}$};

		\end{tikzpicture}
		\caption{The parent square $Q_{k,j}$: five vertical strips $H_{j,1},\ldots,H_{j,5}$, each subdivided into 4 columns. }
		\label{fig:final}
	\end{figure}
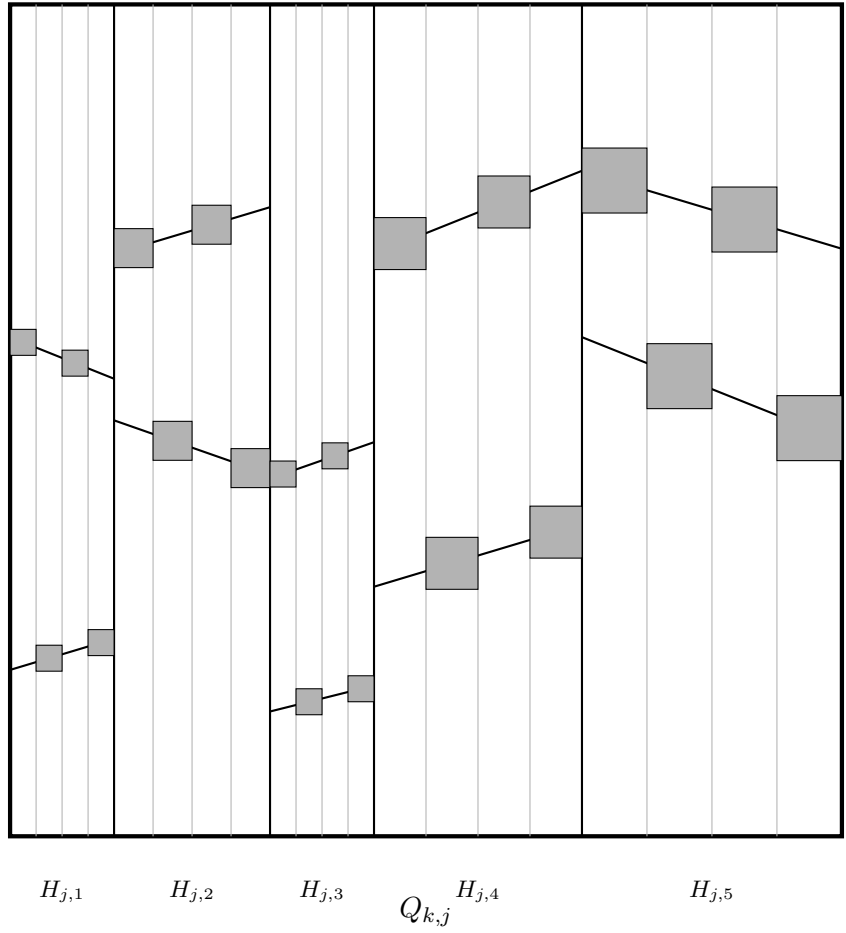

    See Figure \ref{fig:final} for an illustration of our construction.

    We need to check the property claimed to be preserved on each $R_{j,m}^{(i)}=I_{j,m}^{(i)} \times J_{j,m}^{(i)}$. That is, for every $b\in I_{j,m}^{(i)}$, there exists some $(x,y)\in A$ such that $x-yb\in \operatorname{Int}(J_{j,m}^{(i)})$. In our case, one can fix $$(x,y)=\begin{cases}
        (x_{j,m}^+, y_{j,m}^+), & \text{ if $i$ is odd}\\(x_{j,m}^-, y_{j,m}^-), & \text{ if $i$ is even}
    \end{cases}.$$ 
    More precisely, by $|y_{j,m}^\pm|\leq \frac{1}{2}$ and the definition \ref{eq-w}, one can see that the distance between $x_{j,m}^\pm-y_{j,m}^\pm b$ and the center of $J_{j,m}^{(i)}$ is
    $$|w_{(x_{j,m}^\pm, y_{j,m}^\pm)} (b)-w_{(x_{j,m}^\pm, y_{j,m}^\pm)} (b_{j,m}^{(i)})|=y_{j,m}^\pm|b-b_{j,m}^{(i)}|\leq \frac{l_{k+1, j, m}}{4}<\frac{l_{k+1, j, m}}{2}, \ \forall\,b\in I_{j,m}^{(i)}.$$
    
    Then we define $D_{k+1}$ as the union of all child squares $R_{j,m}^{(i)}$ over all $i, j, m$. Notice that the vertical projection $\pi_\infty(D_{k+1})=[0,1]$ remains true because a square is selected in every column.
    
    Finally set $D := \bigcap_{k=0}^\infty D_k$ as a decreasing limit, which is a nonempty compact set.

	We proceed to verify that $D$ satisfies properties 1 and 2. For the former property, because the vertical projection $\pi_\infty(D_k)=[0,1]$ for each $k$, for every $b\in[0,1]$ there exists a cube $Q_{k,j}=I_{k,j}\times J_{k,j}$ in $D_k$ such that $b\in I_{k,j}$. Then, by our construction, there exists $(x_k, y_k)\in A$ such that $(b, x_k-y_kb)\in D_k$. As $A$ and $D$ are compact, every limit point $(x,y)$ of $\{(x_k, y_k)\}_k$ lies in $A$ and satisfies $(b, x-yb)\in D$.
    
    It remains to show that $D$ is an irregular $1$-set. The proof for $0<\mathcal{H}^1(D)<\infty$ is the same as the previous section. We omit details.
    
	For any $p \in D$ and level $k$, choose a radius $r_k = 0.49 \Delta_k$. Let $\mathcal{E}$ be the set of indices $(j,m)$ such that the interval $H_{j,m}$ intersects the vertical projection $\pi_\infty(B(p, r_k))$. Since $2r_k = 0.98 \Delta_k < \Delta_k$, on each $H_{j,m}$, the ball can intersect only one parity family of child squares. 
    
    Let $W_{j,m}$ denote the length of the closed interval $H_{j,m} \cap \pi_\infty(B(p, r_k))$. As intervals $\{H_{j,m}\}_{j,m}$ are almost disjoint, we have $\sum_{(j,m) \in \mathcal{E}} W_{j,m} \le 2r_k$. Also $\pi_\infty(B(p, r_k))$ intersects at  most $\frac{W_{j,m}}{l_{k+1,j,m}} + 2$ many child intervals of $H_{j,m}$ of length $l_{k+1, j, m}$. Due to alternating parity, $B(p, r_k)$ intersects at most $\frac{1}{2}\left(\frac{W_{j,m}}{l_{k+1,j,m}}\right) + 1$ selected child $l_{k+1,j,m}$-squares in $D_{k+1}\cap (H_{j,m}\times \R)$. Hence
		\begin{equation}\begin{aligned}
		    \mathcal{H}^1(D \cap B(p, r_k)) \le & \sum_{(j,m) \in \mathcal{E}} \left( \frac{W_{j,m}}{2 l_{k+1,j,m}} + 1 \right) \sqrt{2}l_{k+1,j,m} \\=& \frac{\sqrt{2}}{2} \sum_{(j,m) \in \mathcal{E}} W_{j,m} + \sqrt{2} \sum_{(j,m) \in \mathcal{E}} l_{k+1,j,m}.
		\end{aligned}
			\end{equation}
			Using $\sum W_{j,m} \le 2r_k$ and bounding the second term by the total number of intervals $H_{j,m}$, we have
			\begin{equation}
				\mathcal{H}^1(D \cap B(p, r_k)) \le \frac{\sqrt{2}}{2}(2r_k) + \sqrt{2} \left( \sum_{j=1}^{M_k} S_j \right) l_{k+1}^{\max}.
			\end{equation}
			Dividing by $2r_k$ and applying \eqref{eq:scale2}:
			\begin{equation}
				\frac{\mathcal{H}^1(D \cap B(p, r_k))}{2r_k} \le \frac{\sqrt{2}}{2} + \sqrt{2} \left( \frac{r_k / 100}{2r_k} \right) \approx 0.707 + 0.007 < 0.72,
			\end{equation}
which implies that $D$ is $1$-irregular.	

	\bibliographystyle{abbrv}
	\bibliography{mybibtex.bib}
	
\end{document}